\documentclass{amsart}

\usepackage{amssymb}
\usepackage{amsmath}
\usepackage{hyperref}

\newcommand{\Gam}{\Gamma_{\mathbb{R}}}
\newcommand{\half}{\frac{1}{2}}
\newcommand{\ep}{\epsilon}
\newcommand{\thalf}{\tfrac{1}{2}}

\newcommand{\sums}{\mathop{{\sum}^{*}}}

\newcommand{\sym}{{\rm sym}}

\numberwithin{equation}{section}

\newtheorem{theorem}{Theorem}[section]

\newtheorem{lemma}[theorem]{Lemma}

\makeatletter
\@namedef{subjclassname@2010}{%
  \textup{2010} Mathematics Subject Classification}
\makeatother

\begin{document}

\title{On the Random Wave Conjecture for Eisenstein series}

\author{Goran Djankovi{\' c}}

\address{
University of Belgrade\\ 
Faculty of Mathematics\\ 
Studentski Trg 16, p.p. 550\\ 
11000 Belgrade, Serbia
}

\email{djankovic@matf.bg.ac.rs}

\author{Rizwanur Khan}

\address{
Department of Mathematics \\ University of Mississippi \\ University, MS 38677-1848, USA}
\email{rrkhan@olemiss.edu}


\begin{abstract} We obtain an asymptotic for the regularized fourth moment of the Eisenstein series for the full modular group, in agreement with the Random Wave Conjecture.
 \end{abstract}
 
 \thanks{The first author was  partially supported by  Ministry of Education,
Science and Technological Development of  Republic of Serbia, Project no.
174008.}

\maketitle

\section{Introduction}

A very interesting and deep conjecture in arithmetic quantum chaos is the Random Wave Conjecture (RWC) of Berry \cite{ber} and Hejhal and Rackner \cite{hejrac}. For the full modular group, the conjecture says that in some sense, the Hecke Maass cusp forms should behave like random waves, in the limit of large Laplacian eigenvalue. This can be formulated more precisely, by proposing that on fixed compact regular sets, the moments of Hecke Maass cusp forms of large eigenvalue, suitably normalized, should asymptotically equal the moments of a standard normal random variable. For the second moment, the statement of the RWC is essentially equivalent to the Quantum Unique Ergodicity (QUE) conjecture, which was resolved quite recently by the work of Lindenstrauss \cite{lin} and Soundararajan \cite{sou} (for other settings, see for example \cite{holsou,nel}). It is notable that the QUE conjecture was actually first settled in the simpler case of the Eisenstein series, by Luo and Sarnak \cite{luoSar}. 

We are interested in the implications of the RWC for the fourth moment, especially because of its relationship with $L$-functions via Watson's formula \cite{wat}. Assuming the Generalized Lindel\"{o}f Hypothesis, Buttcane and the second author \cite{butkha} asymptotically evaluated the fourth moment of Hecke Maass cusp forms over the entire modular curve, and found the result consistent with the RWC. Unconditional works on the fourth moment have only succeeded in obtaining average information (see \cite{blo, kha2}) or upper bounds (see \cite{butkha3, blokhayou, luo, spi}), in various aspects. Obtaining an asymptotic for the fourth moment unconditionally seems out of reach, but would be very important and meaningful in the context of Arithmetic Quantum Chaos.  

In this paper we will focus not on cusp forms but on Eisenstein series. This case seems to be more promising,  but it comes with its own set of difficulties, the first being how to even formulate the problem. Let $E(z,s)$ denote the standard Eisenstein series of weight 0 for $\Gamma=\rm{SL}(2,\mathbb{Z})$.  Let $X=\Gamma\backslash \mathbb{H}$ be the quotient with the upper half complex plane. Given any compact regular subset $\Omega\subset X$, let
\begin{align*}
\tilde{E}(z,\thalf+iT) = \frac{E(z,\thalf+iT)}{ \sqrt{\frac{1}{\text{vol}(\Omega)}\int_{\Omega} |E(z,\thalf+iT)|^2 \ \frac{dxdy}{y^2}}}
\end{align*}
so that $\tilde{E}(z,\thalf+iT)$ has the following normalization on $\Omega$:
\begin{align*}
\frac{1}{\text{vol}(\Omega)}\int_{\Omega} |\tilde{E}(z,\thalf+iT)|^2 \ \frac{dxdy}{y^2}=1.
\end{align*}
The RWC in this case (see \cite[page 300]{hejrac}) is
\begin{align}
\label{compactconje} \frac{1}{\text{vol}(\Omega)}\int_{\Omega} |\tilde{E}(z,\thalf+iT)|^4 \ \frac{dxdy}{y^2} \sim 3,
\end{align}
as $T\to \infty$. The right hand side is the fourth moment of a standard normal random variable. Proving (\ref{compactconje}) for arbitrary compact sets $\Omega$ would probably be out of reach of current methods. One would like to replace $\Omega$ by $X$, but then the integral is not convergent. To get around this difficulty, Spinu \cite{spi} studied the fourth moment of the {\it truncated} Eisenstein series, which decays at the cusp, and he was able to sketch the nice upper bound $O(1)$ for the left hand side (\ref{compactconje}). Some important missing ingredients in Spinu's sketch were provided recently by Humphries \cite{hum}, as we describe below.

We argued in \cite{djakha} that a truncated fourth moment is not the most natural set up. Using Zagier's regularized inner product we instead worked with a regularized fourth moment and formulated a corresponding RWC conjecture for it, the proof of which is our main theorem in this paper. This appears to be the first time that the RWC is verified unconditionally for a fourth moment of an automorphic form.

\begin{theorem}
As $T \rightarrow \infty$, we have 
\begin{align}
\label{oldthmline} \int_{X}^{reg} |E(z , \tfrac{1}{2} + iT)|^4 d \mu(z) \sim \frac{72}{\pi} \log^2 T.
\end{align}
\end{theorem}
\noindent The advantage of the regularized fourth moment is that it has a precise relationship with $L$-functions, which (apart from upper bounds) we did not see in \cite{spi}. This relationship is given in \cite[Theorem 1.1]{djakha}, and using it our task is already reduced to $L$-function theory. Thus to complete the proof of the main theorem, all we need to do is prove

\begin{theorem}\label{propo}
Let $\{ u_j : j\ge 1 \}$ denote an orthonormal basis of even and odd Hecke Maass cusp forms for the full modular group, ordered by Laplacian eigenvalue $\frac{1}{4}+t_j^2$, and let $\Lambda(s, u_j)$ denote the corresponding completed $L$-functions. Let $\xi(s)$ denote the completed Riemann $\zeta$ function. As $T \rightarrow \infty$, we have
\begin{align}
\label{propline} \sum_{\substack{j\ge 1}} \frac{\cosh(\pi t_j)}{2}  \frac{ |\Lambda(\frac{1}{2} + 2Ti, u_j)|^2   \Lambda^2(\frac{1}{2}, u_j) }{L(1, \sym^2 u_j) \, |\xi(1+2Ti)|^4  } \sim \frac{48}{\pi}\log^2 T.
\end{align}
\end{theorem}
\noindent Although we do not display it, it will be clear from the proof that the implied error in (\ref{propline}), and hence also in (\ref{oldthmline}) by \cite[Theorem 1.1]{djakha}, is $O(\log^{5/3+\ep})$.

As we will see, the sum in (\ref{propline}) is essentially supported on $|t_j|<2T+T^\epsilon$. We can analyze the sum over three separate ranges: the short ranges consisting of $|t_j|<T^{1-\epsilon}$ or $||t_j|-2T|<T^{1-\epsilon}$, and the bulk range consisting of $|t_j| \asymp 2T$ but not too close to $2T$. More precisely, the bulk range is
\begin{align}
\label{bulkdef} T^{1-\ep}<|t_j|<2T-T^{1-\ep}.
\end{align}
 While the short ranges were too difficult to treat unconditionally for the fourth moment of Hecke Maass cusp forms, for Eisenstein series they can be quite readily treated using existing work in the literature of Ivi\'{c} \cite{ivi}, Jutila \cite{jut3}, and Jutila and Motohashi \cite{jutmot}. This was observed in \cite{hum} (while a similar treatment for dihedral Maass forms was pointed out in \cite{butkha2}, for one of the short ranges). 

The mean value (\ref{propline}) restricted to the bulk range (\ref{bulkdef}) is reminiscent of the fourth moment of the central values $L(\half, u_j)$ over a dyadic interval of the spectral parameter $t_j$, a problem which has been extensively studied. Such a fourth moment has been generalized by Jutila \cite{jut} to a fourth moment over a shorter interval of $t_j$, and by the second  author \cite{kha3} to the fifth moment (the paper is actually in the weight aspect of holomorphic forms, but the same proof works for the spectral aspect). Our Theorem \ref{propo} can also be viewed as a type of generalization of the fourth moment of $L(\half, u_j)$, more in line with Jutila and Motohashi's \cite{jutmot} hybrid estimates for the mean values of $|L(\half+2Ti, u_j)|^4$ and $L(\half,u_j)^2 | L(\half+2Ti, u_j)|^2$. These results are quite powerful and crucial to the treatment of the short ranges mentioned above. However for the bulk range, the estimates are trivial (for they follow directly from the spectral large sieve), and we must go further to obtain an asymptotic. 

Our treatment of the bulk range is very different from Spinu's sketch of an upper bound.  Most importantly, unlike Spinu (see \cite[section 6.5]{spi}), we completely avoid solving a shifted divisor problem. This is achieved by following the general strategy in \cite{butkha2} and \cite{butkha} of using Kuznetsov's trace formula, Voronoi summation, Kuznetsov's formula again, and finally subconvexity. However the situation here is a bit more delicate because our ``off-diagonal'' contains a main term which must be carefully extracted, as is typical in the study of the fourth moment of $L(\half, u_j)$ and related problems. In fact, and this is quite non-traditional, in our approach the entire main term arises from the off-diagonal. This is because to streamline the analysis we will artificially insert into (\ref{propline}) the root number from the functional equation of $L(\half,u_j)$. This will lead to only the opposite sign Kuznetsov formula, which is the reason why our ``diagonal'' is empty, and because of which we will completely avoid dealing with the $J$-Bessel function and the stationary phase analysis that comes with it (see \cite[section 5.4.9] {spi}). In the context of the fourth moment problem, this trick with the root number was used in a slightly different way but to the same effect in \cite{butkha} and \cite{butkha2}. The trick has been known for some time (see \cite{mot}), and has also appeared more recently in \cite{blolimil, andkir, blo-kha}. 

Throughout we will use the convention that $\ep$ denotes an arbitrarily small positive constant, but not necessarily the same one from one occurrence to the next.

\section{Preliminaries}

\subsection{Reduction to the bulk range} Let \begin{align*}
H(t) =  \frac{ |\Gamma(\frac{\half+2iT+it}{2})|^2 |\Gamma(\frac{\half+2iT-it}{2})|^2 |\Gamma(\frac{\half+it}{2})|^4}{|\Gamma(\frac{1+2iT}{2})|^4  |\Gamma(\frac{1+2it}{2})|^2 }.
\end{align*}
Note that this is the same function from \cite[line (2.3)]{butkha}. By \cite[eqs. (2.5-2.6)]{butkha}, we have that $H(t)$ is essentially supported on $|t|<2T+T^\epsilon$ and is roughly of size $|t|^{-1}|4T^2-t^2|^{-\half}$. Writing out the completed $L$-functions and using the identity $\cosh(\pi T)=\pi /|\Gamma(\frac{1}{2}+iT)|^2$,
we can re-write the left hand side of (\ref{propline}) as 
\begin{align*}
\frac{\pi}{2|\zeta(1+2iT)|^4} \sum_{j\ge 1} H(t_j) \frac{L(\thalf, u_j)^2 |L(\thalf+2iT, u_j)|^2}{{L(1,\sym^2 u_j)}}.
\end{align*}
By \cite[sections 3.6-3.7]{hum}, we have that
\begin{align*}
\frac{\pi}{2|\zeta(1+2iT)|^4}  \Big(\sum_{|t_j|<T^{1-\ep}} + \sum_{2T-T^{1-\ep}<|t|<2T+T^\ep} \Big) H(t_j) \frac{L(\thalf, u_j)^2 |L(\thalf+2iT, u_j)|^2}{{L(1,\sym^2 u_j)}} \ll T^{-\ep}. 
\end{align*}
This uses results of Ivi\'{c} \cite{ivi}, Jutila \cite{jut3}, and Jutila and Motohashi \cite{jutmot} on subconvexity bounds and (hybrid) bounds for moments of Hecke Maass $L$-functions.
Thus it suffices to consider only the bulk range, which we can pick out with the smooth function $W(t)=W_\ep(t)$ given in \cite[Lemma 5.1]{butkha}:
\begin{align}
\label{bulksum} \frac{\pi}{2|\zeta(1+2iT)|^4}  \sum_{j\ge 1} W(t_j) H(t_j) \frac{L(\thalf, u_j)^2 |L(\thalf+2iT, u_j)|^2}{{L(1,\sym^2 u_j)}}.
\end{align}
We recall the main properties of these weight functions:  $H(t)W(t)\ll T^{-100}$ unless $T^{1-\ep}<|t|< 2T-T^{1-\ep}$, in which range
\begin{align*}
\frac{d^k}{dt^k} H(t)W(t) \ll T^{-2} (T^{-1+\ep})^k,
\end{align*}
for all $k \ge 0$
and $W(t)=1+O(T^{-100})$ in the smaller range $T^{1-\ep/4}<|t|< 2T-T^{1-\ep/4}$.  

Let $\lambda_j(n)$ denote the (real) eigenvalues of the $n$-th Hecke operator corresponding to $u_j$, where we write $\lambda_j(-n)=\lambda_j(n)$ for $u_j$ even and $\lambda_j(-n)=-\lambda_j(n)$ for $u_j$ odd. We have that (\ref{bulksum}) equals
\begin{align}
\label{bulksum2} \frac{\pi}{2|\zeta(1+2iT)|^4}  \sum_{j\ge 1} W(t_j) H(t_j) \lambda_j(-1) \frac{L(\thalf, u_j)^2 |L(\thalf+2iT, u_j)|^2}{{L(1,\sym^2 u_j)}}.
\end{align}
This is because $\lambda_j(-1)=1$ if $u_j$ is even, and if $u_j$ is odd then $L(\thalf,u_j)=0$. The purpose of artificially inserting $\lambda_j(-1)$ is to avoid the same sign Kuznetsov formula.
We now work on supplying approximate functional equations for the $L$-values in the sum (\ref{bulksum2}).

\subsection{Approximate functional equations}

The $L$-function attached to $u_j$ is defined for $\Re(s)>1$ by
 \begin{align*}
L(s, u_j) = \sum_{n\ge 1} \frac{\lambda_f(n)}{n^s},
  \end{align*}
with analytic continuation to an entire function. Let $\Gam(s)= \pi^{-\frac{s}{2}} \Gamma(\frac{s}{2})$. For $u_j$ even we have the functional equation (see \cite[chapter 3]{gol})
\begin{align*}
  &L(s, u_j) \Gam(s+it_j)\Gam(s-it_j)= L(1-s, u_j) \Gam(1-s+it_j)\Gam(1-s-it_j),
 \end{align*}
For $u_j$ odd we have
\begin{align*}
L(s, u_j) \Gam(1+s+it_j)\Gam(1+s-it_j)= -L(1-s, u_j) \Gam(2-s+it_j)\Gam(2-s-it_j).
\end{align*}
We have the following standard approximate functional equations (see \cite[Thm 5.3]{iwakow}). For $u_j$ even we have 
\begin{align}
\label{afe1} L(\thalf,u_j)^2 = 2\sum_{nm\ge 1} \frac{\lambda_j(n)\lambda_j(m)}{(nm)^\half}V_1^\text{even}(nm,t_j) = 2\sum_{n, k \ge 1} \frac{\lambda_j(n) \tau(n)}{ k n^\half}V_1^\text{even}(k^2 n,t_j),
\end{align}
where the weight function is defined below, $\tau(n)$ is the divisor function, and the second equality follows from the Hecke relations $\lambda_j(n)\lambda_j(m) = \sum_{k|(n,m)} \lambda_j(nm/k^2)$.
For $u_j$ odd we have
\begin{align}
\label{afe2} L(\thalf,u_j)^2 = 2\sum_{nm\ge 1} \frac{\lambda_j(n)\lambda_j(m)}{(nm)^\half}V_1^\text{odd}(nm,t_j) = 2\sum_{n, k \ge 1} \frac{\lambda_j(n)}{ k n^\half}V_1^\text{odd}(k^2 n,t_j).
\end{align}
For $x, \sigma , |t| >0$, we define
\begin{align*}
&V_1^\text{even}(x,t) = \frac{1}{2\pi i} \int_{(\sigma)} e^{s^2} x^{-s}  \left( \frac{ \Gam(\half +s + it )\Gam(\half +s - it )}{  \Gam(\half  + it )  \Gam(\half - it ) } \right)^2 \frac{ds}{s},\\
&V_1^\text{odd}(x,t) = \frac{1}{2\pi i} \int_{(\sigma)} e^{s^2} x^{-s}  \left( \frac{ \Gam(\frac{3}{2} +s + it )\Gam(\frac{3}{2} +s - it )}{  \Gam(\frac{3}{2}  + it )  \Gam(\frac{3}{2} - it ) } \right)^2 \frac{ds}{s}.
\end{align*}
Using Stirling's approximation of the Gamma function after restricting the integral to $|\Im(s)|<|t|^{\epsilon}$ by the decay of $e^{s^2}$, we see that
\begin{align}
\label{stir} V_1^\text{error}(x,t) := V_1^\text{even}(x,t) - V_1^\text{odd}(x,t)  \ll  |t|^{-1}.
\end{align}
Henceforth, we write
\begin{align*}
V_1(x,t)=V_1^\text{even}(x,t).
\end{align*}
Similarly we have for $u_j$ even that
\begin{align}
\label{afe3} |L(\thalf+2iT, u_j)|^2 = 2\sum_{nm\ge 1} \frac{\lambda_j(n)\lambda_j(m)}{n^{\half+2iT}m^{\half-2iT}}V_2(nm,t_j) = 2\sum_{m, l \ge 1} \frac{\lambda_j(m) \tau(m,2T)}{ l m^\half}V_2(l^2 m,t_j),
\end{align}
where
\begin{align}
\label{taudef} \tau(m,2T) = \sum_{ab=m} \Big(\frac{a}{b}\Big)^{2iT} = \frac{\sigma_{4iT}(m)}{m^{2iT}}, \ \ \ \sigma_k(m)= \sum_{d|m} d^k,
\end{align}
and 
\begin{align*}
V_2(x,t) = \frac{1}{2\pi i} \int_{(\sigma)} e^{s^2} x^{-s}  \prod_\pm \frac{ \Gam(\half + s +i(\pm 2T+t)   )\Gam(\half +  s +i(\pm 2T- t) )}{  \Gam(\half+  i(\pm 2T + t) )  \Gam(\half + i(\pm 2T  - t) ) }  \frac{ds}{s}.
\end{align*}
By Stirling's estimates, we have
\begin{align}
\label{stir2} V_1(x,t)\ll_\sigma \Big(\frac{|t|^2+1}{x}\Big)^{\sigma}, \ \ \ V_2(x,t)\ll_\sigma \Big(\frac{|4T^2-t^2|+1}{x}\Big)^{\sigma}
\end{align}
for any $\sigma>0$.

\bigskip

We are ready to apply the approximate functional equations. By (\ref{afe3}), we have that (\ref{bulksum2}) equals
\begin{align*}
  \frac{\pi}{|\zeta(1+2iT)|^4}  \sum_{j\ge 1} \frac{W(t_j) H(t_j)}{L(1,\sym^2 u_j)}  \lambda_j(-1) L(\thalf, u_j)^2 \sum_{m, l \ge 1} \frac{\lambda_j(m) \tau(m,2T)}{ l m^\half}V_2(l^2 m,t_j).
\end{align*}
We were able to apply (\ref{afe3}) even though it holds only for $u_j$ even because otherwise $L(\thalf, u_j)=0$.  Now using (\ref{afe1}-\ref{afe2}) we that have this equals
\begin{multline}
 \label{bulksum3} \frac{2\pi}{|\zeta(1+2iT)|^4}  \sum_{j\ge 1} \frac{W(t_j) H(t_j)}{L(1,\sym^2 u_j)}  \sum_{n,m, k,l \ge 1}  \frac{\lambda_j(-1) \lambda_j(n) \lambda_j(m) \tau(n) \tau(m,2T)}{ k l (nm)^\half} V_1(k^2 n, t_j) V_2(l^2 m,t_j) \\+ \text{Error},
\end{multline}
where Error denotes a similar sum, but with $V_1(k^2 n, t_j)$ replaced by $V_1^\text{error}(k^2 n, t_j)$. By (\ref{stir2}), both the $n$ and $m$ sums have length essentially $T^{2+\ep}$ (remember that $t_j$ is restricted to the bulk range), so by the spectral large sieve \cite[Theorem 7.24]{iwakow}, we immediately get the upper bound $O(T^\ep)$ for the first sum of (\ref{bulksum3}). By (\ref{stir}) and the spectral large sieve we get that Error is $O(T^{-1+\ep})$.  Thus the main result (\ref{propline}) is reduced to proving that
\begin{align}
\label{reduced} \sum_{j\ge 1} \frac{W(t_j) H(t_j)}{L(1,\sym^2 u_j)}  \sum_{n,m, k,l \ge 1}  \frac{ \lambda_j(-n) \lambda_j(m) \tau(n) \tau(m,2T)}{ k l (nm)^\half} V_1(k^2 n, t_j) V_2(l^2 m,t_j) \sim \frac{24}{\pi^2} |\zeta(1+2iT)|^4 \log^2 T.
\end{align}
Let us keep in mind (see \cite[Theorem 8.27, 8.29]{iwakow}) that $\zeta(1+2iT)$ and $1/\zeta(1+2iT)$ are bounded above by some power of $\log T$.

\subsection{Kuznetsov's trace formula}

Applying Kuznetsov's trace formula (we follow the normalization given in \cite[Lemma 3.2]{butkha}), we get that the left hand side of (\ref{reduced}) equals $\mathcal{E}+\mathcal{O}$, where $\mathcal{E}$ is the Eisenstein series contribution and $\mathcal{O}$ is the off-diagonal contribution, as given below. Note that there is no diagonal contribution because of the opposite signs in $\lambda_j(-n)$ and $\lambda_j(m)$: thus the main term must arise from elsewhere.

By the same type of argument as \cite[section 5]{butkha2}, we have that the Eisenstein series contribution is bounded as follows:
\begin{align*}
\mathcal{E}\ll  \int_{-\infty}^{\infty} H(t)W(t) \frac{|\zeta(\half+it)|^4  |\zeta(\half+i(2T-t))|^4  }{|\zeta(1+2it)|^2} dt.
\end{align*}
The weight function $W(t)$ restricts the integral to the bulk range, in which we have already noted that $W(t)H(t) \ll T^{-2+\ep}$. Using any subconvexity bound for the four factors of $\zeta(\half+i(2T-t))$, and the large sieve \cite[Theorem 9.1]{iwakow} to sharply bound the fourth moment of $\zeta(\half+it)$, we get that $\mathcal{E}\ll T^{-\delta}$ for some $\delta>0$.

It remains to consider
\begin{multline}
\label{rem} \mathcal{O}= \frac{2}{\pi^2}  \sum_{n,m, k,l, c \ge 1}   \frac{ \tau(n) \tau(m,2T)}{ k l (nm)^\half} \frac{S(-n,m,c)}{c} \\ \int_{-\infty}^{\infty} \sinh(\pi t) K_{2it}\Big(\frac{4\pi\sqrt{nm}}{c}\Big) V_1(k^2 n, t) V_2(l^2 m,t)H(t)W(t) \ tdt,
\end{multline}
and we need to show that
\begin{align}
 \mathcal{O}\sim  \frac{24}{\pi^2} |\zeta(1+2iT)|^4 \log^2 T.
\end{align}
Before we evaluate the integral transform, we make some simplifications. 

By the rapid decay of $W(t)$, we may insert a function $Z(\frac{t}{2T})$ in the integrand above, where $Z(x)$ is even, compactly supported on $T^{-2\ep}<|x|< 1-T^{-2\ep}$ with 
\begin{align}
\label{zprop} Z(x)=1 \  \text{ for } \ T^{-\ep}<|x|< 1-T^{-\ep},
\end{align}
 and smooth with with derivatives $\|Z^{(r)}\|_\infty\ll (T^{\ep})^r$. Henceforth, we absorb $W(t)$ into $Z(\frac{t}{2T})$. Note that we could not make this simplification earlier because we neeeded a carefully constructed weight function such as $W(t)$ which is admissible in Kuznetsov's formula.  By \cite[lines (2.5) and (5.4)]{butkha}, we may replace $H(t)$ by the leading term in its Stirling's expansion as the lower order terms can be treated similarly:
\begin{align*}
H(t)= \frac{8\pi}{|t|(4T^2-t^2)^\half} + \ldots.
\end{align*}
Also by Stirling's approximation (similarly to \cite[lines (3.12-3.13)]{butkha}), we may treat only the leading terms of $V_1(k^2 n, t)$ and  $V_2(l^2 m,t)$:
\begin{align*}
V_1(k^2 n, t)= V\Big(\frac{k^2n}{t^2}\Big)+\ldots, \ \ \ V_2(l^2 m, t)= V\Big(\frac{l^2 m}{4T^2-t^2}\Big)+\ldots,
\end{align*}
where
\begin{align}
\label{vdef} V(x)=\frac{1}{2\pi i} \int_{(c)} e^{s^2} (4\pi^2 x)^{-s} \frac{ds}{s} +\ldots
\end{align}
Now by \cite[Lemma 3.4]{butkha}, we can evaluate the integral transform (treating the leading order term) in (\ref{rem}) as follows:
\begin{multline}
\label{between2} \int_{-\infty}^{\infty} \sinh(\pi t) K_{2it}\Big(\frac{4\pi\sqrt{nm}}{c}\Big) Z\Big(\frac{t}{2T} \Big) V\Big(\frac{k^2n}{t^2}\Big) V\Big(\frac{l^2 m}{4T^2-t^2}\Big) \frac{8\pi}{|t|(4T^2-t^2)^\half}  \ tdt  \\
=\frac{8\pi^3 \sqrt{nm}}{c} Q\Big(\frac{2\pi \sqrt{nm} }{c}\Big) V\Bigg(\frac{ k^2 c^2}{4\pi^2 m}\Bigg) V\Bigg(\frac{ l^2 m }{4T^2} \frac{1}{(1 - \frac{\pi^2 nm}{T^2 c^2})}\Bigg) + \ldots,
\end{multline}
where
\begin{align}
\label{qdef} Q(t)=  \frac{Z(\frac{t}{2T})}{|t|(4T^2-t^2)^\half}
\end{align}
and the ellipsis denotes lower order terms and a negligible error incurred from applying the lemma (this error is a negative power of $T$ by the same argument as \cite[line (8.2)]{butkha2}). Thus by (\ref{rem}-\ref{between2}), we need to prove
\begin{multline}
\label{rem1} \sum_{n,m, k,l, c \ge 1}   \frac{ \tau(n) \tau(m,2T) S(-n,m,c)}{ k l c^2} Q\Big(\frac{2\pi \sqrt{nm} }{c}\Big)  V\Bigg(\frac{ k^2 c^2}{4\pi^2 m}\Bigg) V\Bigg(\frac{ l^2 m }{4T^2} \frac{1}{(1 - \frac{\pi^2 nm}{T^2 c^2})}\Bigg) \\ 
 \sim  \frac{3}{2\pi^3} |\zeta(1+2iT)|^4 \log^2 T.
\end{multline}
In this sum, we can assume by the properties of $V$ and $Q$ that
\begin{align}
\label{ranges} m\le \frac{T^{2+\ep}}{l^2}, \ \ \  \frac{T^{2-\ep} c^2}{m} \le n \le \frac{T^{2+\ep} c^2}{m} \le \frac{T^{2+\ep} }{k^2}, \ \ \ \frac{\sqrt{nm}}{T^{1+\ep}} \le c \le \frac{\sqrt{nm}}{T^{1-\ep}} \le T^{2+\ep}.
\end{align}

\subsection{Voronoi summation} We will need the Voronoi summation formula:

 \begin{lemma} \cite[Lemma 3.1]{kha3} \label{vorsum1} 
Given a smooth function $\Phi$, compactly supported on the positive reals, and coprime integers $h$ and $c$, we have
\begin{align*}
\sum_{n\ge 1} \frac{\tau(n)}{n} e\Big(\frac{n\overline{h}}{c}\Big)\Phi\Big(\frac{n}{N}\Big) = \mathcal{M} + \mathcal{D},
\end{align*}
where the main term $\mathcal{M}$ and dual sum $\mathcal{D}$ are given by
\begin{align*}
\mathcal{M}= \frac{1}{c} \int_{-\infty}^\infty \Big(\log\frac{ x}{c^2}+2\gamma\Big)\Phi \Big(\frac{x}{N}\Big)\frac{dx}{x}, \ \ \ \ \ \ \  \mathcal{D}= \sum_\pm \frac{1}{c} \sum_{r\ge 1} \tau(r) e\Big(\frac{\pm rh}{c}\Big) \check{\Phi}_\pm \Big(\frac{Nr}{c^2}\Big),
\end{align*}
where
\begin{align}
\label{phi-trans}\check{\Phi}_\pm(x)= \frac{1}{2\pi i} \int_{(\sigma)} G_1^\pm(s) \tilde{\Phi}(-s)  x^{-s}  ds, \ \ G_1^\pm(s)= 2(2\pi )^{-2s} \Gamma(s)^2  \cos^{(1\mp1)/2} (\pi s),
\end{align}
$\tilde{\Phi}$ is the Mellin transform of $\Phi$, $\sigma>0$, and $\gamma$ is Euler's constant.
\end{lemma}
\noindent By Stirling's approximation, for $\Re(s)=\sigma$ fixed away from $\mathbb{Z}_{\le 0}$, we have 
\begin{align}
\label{g1bound} |G_1^\pm(s)|\ll_\sigma (1+|\Im(s)|)^{2\sigma-1}
\end{align}

We will also need a Voronoi summation formula for $\tau(m,2T)$.

 \begin{lemma} \cite[Lemma 3.3]{hou} \label{vorsum2} 
Given a smooth function $\Phi$, compactly supported on the positive reals, and coprime integers $h$ and $c$, we have
\begin{align*}
\sum_{m\ge 1} \frac{\tau(m,2T)}{m} e\Big(\frac{m h }{c}\Big)\Phi\Big(\frac{m}{M}\Big) = \mathcal{M} + \mathcal{D},
\end{align*}
where the main term $\mathcal{M}$ and dual sum $\mathcal{D}$ are given by
\begin{align}
\label{vor2md} &\mathcal{M}= \frac{\zeta(1+4Ti)}{c^{1+4Ti}} \int_{-\infty}^\infty \Phi(x) x^{2Ti} \frac{dx}{x} + \frac{\zeta(1-4Ti)}{c^{1-4Ti}} \int_{-\infty}^\infty \Phi(x) x^{-2Ti} \frac{dx}{x}, \\ 
\nonumber &\mathcal{D}=   \sum_\pm \frac{1}{c} \sum_{b\ge 1} \tau(b,2T) e\Big(\frac{\pm b\overline{h}}{c}\Big) \check{\Phi}^\pm\Big(\frac{Mb}{c^2}, T\Big),
\end{align}
where
\begin{align}
\label{phi-trans2} &\check{\Phi}^\pm(x,T)= \frac{1}{2\pi i} \int_{(\sigma)} G_2^\pm(s) \tilde{\Phi}(-s)  x^{-s}  ds, \\
&\nonumber G_2^\pm(s)=  2(2\pi )^{-2s} \Gamma(s-2Ti) \Gamma(s+2Ti) \cos^{(1\mp1)/2}(\pi s) \cosh^{(1\pm1)/2}(2\pi T),
\end{align}
$\tilde{\Phi}$ is the Mellin transform of $\Phi$, and $\sigma>0$.
\end{lemma}

%\noindent This formula can be proven in the same way as \cite[Theorem 1.7]{jut2}, using the Estermann Zeta function. In fact the formula appears (apart from a couple of errors which we have corrected) in an unpublished appendix to a paper of Kowalski, Michel, and Vanderkam \cite[page 48]{kmv}.

By Stirling's approximation, for $\Re(s)=\sigma$ fixed away from $\mathbb{Z}_{\le 0}$, we have 
\begin{align}
\label{g2bound} |G_2^\pm(s)|\ll_\sigma (1+|\Im(s)-2T|)^{\sigma-\half} (1+|\Im(s)+2T|)^{\sigma-\half}.
\end{align}
If we also have $\Im(s)\ll T^\ep$, then
\begin{align*}
 G_2^\pm(s) = 2 (T/\pi)^{2s-1} e^{-2\pi T} \cos^{(1\mp1)/2}(\pi s) \cosh^{(1\pm1)/2}(2\pi T)(1+O(T^{-1+\ep}).
\end{align*}
See \cite[section 3.2]{butkha2} for more details on these calculations. Thus only the + case is significant in this range of $s$, and we have
\begin{align}
\label{g2asymp2} G_2^+(s) =  (T/\pi)^{2s-1} (1+O(T^{-1+\ep}).
\end{align}

\section{The main term}

We divide the sum (\ref{rem1}) into dyadic intervals of $n$ by using a partition of unity. Thus we need to consider the sum
\begin{multline}
\label{rem2} \sum_j \sum_{n,m, k,l, c \ge 1}  \Psi_{1,j}\Big(\frac{n}{N_j}\Big)   \frac{ \tau(n) \tau(m,2T) S(-n,m,c)}{ k l c^2} Q\Big(\frac{2\pi \sqrt{nm} }{c}\Big)  V\Bigg(\frac{ k^2 c^2}{4\pi^2 m}\Bigg) V\Bigg(\frac{ l^2 m }{4T^2} \frac{1}{(1 - \frac{\pi^2 nm}{T^2 c^2})}\Bigg)
\end{multline}
for a sequence of smooth functions $\Psi_{1,j}$ compactly supported on $(\half, \frac{3}{2})$ say. We open the Kloosterman sums, writing
\begin{align*}
S(-n,m,c)=\sums_{h\bmod c} e\Big(\frac{nh-m\overline{h}}{c}\Big)
\end{align*}
and then apply Voronoi summation (Lemma \ref{vorsum1}) to the sum over $n$. The main term $\mathcal{M}$ resulting from this lemma  equals, after reversing the partition of unity,
\begin{align*}
\mathcal{M}= \sum_{m, k,l, c \ge 1}   \frac{  \tau(m,2T) r_c(m) }{ k l c^3} V\Bigg(\frac{ k^2 c^2}{4\pi^2 m}\Bigg)  \int_{0}^{\infty} \Big( \log\frac{x}{c^2}+2\gamma\Big)Q\Big(\frac{2\pi \sqrt{xm} }{c}\Big)   V\Bigg(\frac{ l^2 m }{4T^2} \frac{1}{(1 - \frac{\pi^2 xm}{T^2 c^2})}\Bigg)  dx,
\end{align*}
where $r_c(m)=\sum_{h\bmod c} e(m\overline{h}/c)$ is the Ramanujan sum. The goal of this section is to show that 
\begin{align}
\label{mt} \mathcal{M}\sim  \frac{3}{2\pi^3} |\zeta(1+2iT)|^4 \log^2 T,
\end{align}
which is the main term expected in (\ref{rem1}). In the following sections we will show that the dual sum arising from Voronoi summation is $O(T^{-\delta})$ for some $\delta>0$, and this will complete the proof of the main result.

\bigskip

Making the substitution $y=\frac{2\pi \sqrt{xm} }{c}$, we have
\begin{align*}
 \mathcal{M}= \frac{1}{2\pi^2}  \sum_{m, k,l, c \ge 1} \frac{  \tau(m,2T) r_c(m) }{ k l c m} V\Bigg(\frac{ k^2 c^2}{4\pi^2 m}\Bigg)  \int_{0}^{\infty} \Big( \log\frac{y^2}{4\pi^2 m}+2\gamma\Big) Q(y)   V\Bigg(\frac{ l^2 m }{4T^2} \frac{1}{(1 - \frac{y^2}{4T^2})}\Bigg)  y dy.
\end{align*}
This is trivially bounded by $O(T^\ep)$. Now we insert (\ref{qdef}) and up to an error of $O(T^{-\ep})$, we extend $Z(\frac{y}{2T})$ to the endpoints 0 and $2T$, getting
\begin{align*}
\mathcal{M} \sim \frac{1}{2 \pi^2} \sum_{m, k,l, c \ge 1}  \frac{  \tau(m,2T) r_c(m) }{ k l c m} V\Bigg(\frac{ k^2 c^2}{4\pi^2 m}\Bigg)  \int_0^{2T} \Big(\log\frac{y^2}{4\pi^2 m} + 2\gamma \Big) \frac{1}{(4T^2-y^2)^\half}   V\Bigg(\frac{ l^2 m }{4T^2} \frac{1}{(1 - \frac{y^2}{4T^2})}\Bigg)  dy.
\end{align*}
Making the substitution $x=\frac{y}{2T}$, and keeping only the leading contribution from the log, we get
\begin{align*}
\mathcal{M} \sim  \frac{1}{2 \pi^2} \sum_{m, k,l, c \ge 1}  \frac{  \tau(m,2T) r_c(m) }{ k l c m} V\Bigg(\frac{ k^2 c^2}{4\pi^2 m}\Bigg)  \int_0^{1} \log\Big(\frac{T^2}{m}\Big) \frac{1}{(1-x^2)^\half}   V\Bigg(\frac{ l^2 m }{4T^2} \frac{1}{(1 - x^2)}\Bigg)  dx.
\end{align*}
Using (\ref{vdef}), we have
\begin{multline}
\label{rem3} \mathcal{M}\sim  \frac{1}{2\pi^2} \frac{1}{(2\pi i)^2} \int_{(\sigma_1)} \int_{(\sigma_2)} e^{s_1^2+s_2^2}  \pi^{-2s_2} T^{2s_2} \sum_{m, k,l, c \ge 1}  \frac{  \tau(m,2T) r_c(m) }{ \ k^{1+2s_1} l^{1+2s_2} c^{1+2s_1} m^{1-s_1 + s_2}} \\ \int_0^{1} \log\Big(\frac{T^2}{m}\Big) \frac{dx}{(1-x^2)^{\half-s_2}}      \ \frac{ds_1}{s_1} \frac{ds_2}{s_2},
\end{multline}
where for absolute convergence we require $\sigma_1, \sigma_2 >0$ and $\sigma_2>\sigma_1$.
Using (\ref{taudef}) and the well known identities
\begin{align*}
\sum_{c\ge 1} \frac{r_c(m)}{c^s} = \frac{\sigma_{s-1}(m)}{m^{s-1} \zeta(s)}, \ \ \ \ \sum_{m\ge 1} \frac{\sigma_a(m) \sigma_b(m)}{m^s} = \frac{\zeta(s)\zeta(s-a)\zeta(s-b)\zeta(s-a-b)}{\zeta(2s-a-b)},
\end{align*}
we get that
\begin{align*}
 &\sum_{m, k, c \ge 1}  \frac{  \tau(m,2T) r_c(m) }{ k^{1+2s_1}  c^{1+2s_1} m^{1-s_1 + s_2}} = F(s_1,s_2),\\
 &\sum_{m, k, l,c \ge 1}  \frac{  \tau(m,2T) r_c(m) }{ k^{1+2s_1} l^{1+2s_2} c^{1+2s_1} m^{1-s_1 + s_2}} = F(s_1,s_2)\zeta(1+2s_2),
\end{align*}
where
\begin{align*}
F(s_1,s_2) = \frac{\zeta(1+s_1+s_2+2iT)\zeta(1+s_1+s_2-2iT)\zeta(1-s_1+s_2+2iT)\zeta(1-s_1+s_2-2iT) }{\zeta(2+2s_2)  }.
\end{align*}
Now writing $\log(T^2/m)=2\log T -\log m$, we have $\mathcal{M}\sim \mathcal{M}_1+\mathcal{M}_2$, where
\begin{align*}
&\mathcal{M}_1=  \frac{\log T}{\pi^2} \frac{1}{(2\pi i)^2} \int_{(\sigma_1)} \int_{(\sigma_2)} e^{s_1^2+s_2^2}  \pi^{-2s_2} T^{2s_2} F(s_1,s_2) \zeta(1+2s_2)  \int_0^{1}  \frac{dx}{(1-x^2)^{\half-s_2}}      \ \frac{ds_1}{s_1} \frac{ds_2}{s_2},\\
&\mathcal{M}_2 =\frac{1}{2\pi^2} \frac{1}{(2\pi i)^2} \int_{(\sigma_1)} \int_{(\sigma_2)} e^{s_1^2+s_2^2}  \pi^{-2s_2}  T^{2s_2} \Big( \frac{d}{ds_2} F(s_1,s_2) \Big) \zeta(1+2s_2) \int_0^{1}  \frac{dx}{(1-x^2)^{\half-s_2}}      \ \frac{ds_1}{s_1} \frac{ds_2}{s_2}.
\end{align*}
We move the $s_2$-integrals of $\mathcal{M}_1$ and $\mathcal{M}_2$ just left of the line $\Re(s_2)=0$. The integrals on the left are 
bounded by a negative power of $T$, so we need only consider the residues from the poles.

In the case of $\mathcal{M}_1$ we crossed a double pole at $s_2=0$ and two simple poles at $s_2=s_1 \pm 2iT$. The contribution of the simple poles is exponentially small as a function of $T$ because of the rapid decay of $e^{s_2^2}$ for $|\Im(s_2)|\to\infty$. 
So we need only consider  the residue of the double pole at $s_2=0$ in $\mathcal{M}_1$. Writing $\zeta(1+2s_2)=\frac{1}{2s_2}+\gamma+\ldots$, we see that the residue at $s_2=0$ is given by 
\begin{align}
\label{resm1} \text{res}_{\mathcal{M}_1} = &  \frac{  \gamma \log T}{\pi^2 } \frac{1}{2\pi i} \int_{(\sigma_1)} e^{s_1^2} F(s_1, 0 ) \int_0^1 \frac{dx}{(1-x^2)^\half}   \frac{ds_1}{s_1}\\
\nonumber + & \frac{ \log T}{2\pi^2 } \frac{1}{2\pi i} \int_{(\sigma_1)} e^{s_1^2} \int_0^1 \frac{1}{(1-x^2)^\half} \Big(  \frac{d}{ds_2}\Big|_{s_2=0}  \frac{e^{s_2^2}  T^{2s_2} F(s_1,s_2)}{\pi^{2s_2}(1-x^2)^{-s_2} }\Big) \ dx \frac{ds_1}{s_1}. 
\end{align}
In the second line, the contribution from the derivative of $T^{2s_2}$ is
\begin{align}
& \nonumber \frac{  \log^2 T}{\pi^2 } \frac{1}{2\pi i} \int_{(\sigma_1)} e^{s_1^2} F(s_1, 0 ) \int_0^1 \frac{dx}{(1-x^2)^\half}   \frac{ds_1}{s_1}\\
 \label{rem4} &= \frac{  \log^2 T}{ 2\pi } \frac{1}{2\pi i} \int_{(\sigma_1)} e^{s_1^2} F(s_1, 0 )   \frac{ds_1}{s_1}. 
\end{align}
To evaluate the $s_1$-integral, we observe that $e^{s_1^2} F(s_1, 0 )$ is an even function, so moving the line of integration to the left of $\Re(s_1)=0$, picking up a residue at $s_1=0$, and then making the substitution $s_1 \leftrightarrow -s_1$, we get that the $s_1$-integral equals half of the residue of its integrand. Thus (\ref{rem4}) equals
\begin{align*}
 \frac{  \log^2 T}{2\pi } \cdot \frac{F(0,0)}{2} =  \frac{  \log^2 T}{ 2\pi } \cdot \frac{|\zeta(1+2iT)|^4}{2\zeta(2)}  =  \frac{3}{ 2\pi^3}  |\zeta(1+2iT)|^4 \log^2 T.
\end{align*}
This dominates by a factor of $\log T$ the contribution of the first line in (\ref{resm1}) and the contribution of the derivative of $e^{s_2^2} \pi^{-2s_2} (1-x^2)^{-s_2}$ in the second line. It also dominates by a lesser power of $\log T$ the contribution of the derivative of   $F(s_1,s_2)$ in the second line of (\ref{resm1}) because 
\begin{align} \label{zeta_deriv}
\frac{\zeta'(1+2iT)}{\zeta(1+2iT)} \ll (\log T)^{\frac23+\ep}
\end{align}
by classical estimates \cite[Theorem 8.29]{iwakow}. Thus so far we have shown
\begin{align*}
{\mathcal{M}_1} \sim \frac{3}{2\pi^3}  |\zeta(1+2iT)|^4 \log^2 T,
\end{align*}
which is exactly the expected main term in (\ref{mt}). 

Indeed,  we now turn to the contribution of the $s_1$-integrals of residues in ${\mathcal{M}_2}$ and show that it is asymptotically smaller. Here we encountered double poles  at $s_2=s_1 \pm 2iT$ (arising from the derivative of $\zeta(1-s_1+s_2 \pm 2iT)$), whose contribution is exponentially small in $T$ as before, and a double pole at $s_2=0$, the residue of which we denote by $\text{res}_{\mathcal{M}_2}$. We have explicitly 
\begin{align*}
\text{res}_{\mathcal{M}_2} &=\frac{\gamma}{2\pi^2} \frac{1}{2\pi i} \int_{(\sigma_1)}  e^{s_1^2}   \Big( \frac{d}{ds_2}\Big|_{s_2=0} F(s_1,s_2) \Big)  \int_0^{1}  \frac{dx}{(1-x^2)^{\half}}      \ \frac{ds_1}{s_1} \\
& \quad +\frac{1}{4\pi^2 } \frac{1}{2\pi i} \int_{(\sigma_1)} e^{s_1^2} \int_0^1 \frac{1}{(1-x^2)^\half} \Big(   \frac{d}{ds_2}\Big|_{s_2=0}  \frac{e^{s_2^2}  T^{2s_2} \frac{d}{ds_2} F(s_1,s_2)}{\pi^{2s_2}(1-x^2)^{-s_2} }\Big) \ dx \frac{ds_1}{s_1}.
\end{align*}
Again using (\ref{zeta_deriv}), we obtain that the contribution of the first line is $ O( |\zeta(1+2iT)|^4 (\log T)^{\frac23+\ep})$. In the second line, the derivative of $T^{2s_2}$  contributes a factor $2\log T$,   giving contribution of $ O( |\zeta(1+2iT)|^4 (\log T)^{\frac{5}{3}+\ep})$. In a similar way (see \cite[Lemma 4.3]{djakha}), the term with the second derivative $\frac{d^2}{ds_2^2}\Big|_{s_2=0} F(s_1,s_2)$ contributes $O( |\zeta(1+2iT)|^4 (\log T)^{\frac{4}{3}+\ep})$. All other terms in the second line are even smaller.

All in all,
\begin{align*}
\mathcal{M} = \mathcal{M}_1+\mathcal{M}_2  \sim \mathcal{M}_1 \sim \frac{3}{2\pi^3}  |\zeta(1+2iT)|^4 \log^2 T,
\end{align*}
as desired in (\ref{mt}).

\section{The error terms}

\subsection{First application of Voronoi summation: the dual sum} We now return to (\ref{rem2}), to which we applied Voronoi summation in the $n$-sum. The main term arising from the summation formula was already treated, so now we focus on the dual sum. Since there are $O(T^\ep)$ dyadic intervals in (\ref{rem2}) it suffices to consider any one (let $\Psi_{1,j}=\Psi_1$ and $N_j=N$). The dual sum (see Lemma \ref{vorsum1}) is 
\begin{align*}
\mathcal{D} =  N \sum_{\pm}  \sum_{r,m, k,l, c \ge 1} \sums_{h \bmod c} e\Big(\frac{h(-m\pm r)}{c}\Big)   \frac{ \tau(r) \tau(m,2T) }{ k l c^3}  V\Bigg(\frac{ k^2 c^2}{4\pi^2 m}\Bigg)   \check{\Phi}_\pm \Big(\frac{Nr}{c^2}\Big),
\end{align*}
where $\check{\Phi}_\pm$ was defined in (\ref{phi-trans}), 
\begin{align*}
\Phi(x)= x\Psi_1(x)  Q\Big(\frac{2\pi \sqrt{xNm} }{c}\Big) V\Bigg(\frac{ l^2 m }{4T^2} \frac{1}{(1 - \frac{\pi^2 x Nm}{T^2 c^2})}\Bigg),
\end{align*}
$\Psi_1$ is any smooth function compactly supported on $(\frac{1}{2},\frac32)$ with derivatives $\|\Psi_1^{(r)}\|_\infty\ll (T^{\ep})^r$, and $1\le N \le  T^{2+\ep}$.  The goal is to prove $\mathcal{D}\ll T^{-\delta}$ for some $\delta>0$.

\bigskip

We simplify the notation a bit before moving on. By (\ref{qdef}), we have
\begin{align*}
Q(t)=  \frac{1}{4T^2} \frac{Z(\frac{t}{2T})}{\frac{|t|}{2T}(1-(\frac{t}{2T})^2)^\half}.
\end{align*}
We can redefine $Z(\frac{t}{2T})$ and simply write $Q(t)$ as $\frac{1}{T^2}Z(\frac{t}{2T})$. By doing this we lose property (\ref{zprop}) but this is irrelevant for us now as we are aiming for an upper bound. Thus
\begin{align*}
\Phi(x)= \frac{x\Psi_1(x) }{T^2}  Z\Big(\frac{\pi \sqrt{xNm} }{cT}\Big) V\Bigg(\frac{ l^2 m }{4T^2} \frac{1}{(1 - \frac{\pi^2 x Nm}{T^2 c^2})}\Bigg).
\end{align*}
Now referring to (\ref{phi-trans}) for the definition of $ \check{\Phi}_\pm(\frac{Nr}{c^2})$, we observe that in this definition we have 
\begin{align}
\label{phimell} \tilde{\Phi}(-s) = \int_0^\infty \Phi(x) x^{s-1} dx \ll T^\ep \Big(\frac{T^{\ep}}{1+|s|}\Big)^A
\end{align}
 for any $A\ge 0$ by integrating by parts, and by moving the $s$-integral in (\ref{phi-trans}) far to the right (taking $\sigma$ large) and using (\ref{g1bound}) we can assume that $\frac{Nr}{c^2}\ll T^\ep$. Thus it suffices to restrict to $|\Im(s)|<T^\ep$ with $\Re(s)=\ep$ say, and to 
 \begin{align*}
 r< \frac{c^2T^{\ep}}{N} < \frac{c^2T^{\ep}}{c^2 T^2/m} < \frac{m}{T^{2-\ep}}
 \end{align*}
 by (\ref{ranges}). Thus unless $m$ is at its maximum range, we are done (because the $r$-sum would be empty). So henceforth we assume that $1\le r < T^\ep$ and introduce a smooth bump function $\Psi_2(\frac{m}{M})$ to restrict the range of $m$. We have shown that it suffices to prove that
\begin{multline*}
  \frac{N}{T^2} \sum_{\pm}  \sum_{m, k,l, c \ge 1} \sums_{h \bmod c} e\Big(\frac{h(-m\pm r)}{c}\Big)   \frac{  \tau(m,2T) }{ k l c^3} x^{-s} \Psi_1(x)  \Big(\frac{Nr}{c^2}\Big)^{-s} Z\Big(\frac{\pi \sqrt{xNm} }{cT}\Big) \Psi_2\Big(\frac{m}{M}\Big) \\  V\Bigg(\frac{ k^2 c^2}{4\pi^2 m}\Bigg) V\Bigg(\frac{ l^2 m }{4T^2} \frac{1}{(1 - \frac{\pi^2 x Nm}{T^2 c^2})}\Bigg) \ll T^{-\delta}
\end{multline*}
for any given
\begin{align}
\label{ranges0} r<T^\ep, \ \ \ \ T^{2-\ep}<M<T^{2+\ep},
\end{align}
and $s$ with $|\Im(s)|<T^\ep$ and $\Re(s)=\ep$. We can drop $ x^{-s} \Psi_1(x) $ and absorb $(\frac{Nr}{c^2})^{-s} $ into the existing weight functions. By (\ref{ranges}), the function $Z$ restricts $c$ to $N^\half T^{-\ep} < c < N^\half T^{\ep}$. 

We continue to simplify. We can write $V$ as an integral using (\ref{vdef}). In this formula, we can restrict to $\Re(s)=\ep$, which is just right of the integrand's pole at $s=0$, and $|\Im s|<T^\ep$ by the rapid decay of $e^{s^2}$. Thus it suffices to prove
\begin{multline*}
  \frac{N}{T^2} \sum_{\pm}  \sum_{m, k,l, c \ge 1} \sums_{h \bmod c} e\Big(\frac{h(-m\pm r)}{c}\Big)   \frac{  \tau(m,2T) }{ k l c^3}  Z\Big(\frac{\sqrt{Nm} }{cT}\Big)  \Psi_2\Big(\frac{m}{M}\Big)  \\ \Big(\frac{k^2 c^2}{m}\Big)^{-s_1} \Big(\frac{ l^2 m }{T^2}\Big)^{-s_2} \Big(1 - \frac{\pi^2 x Nm}{T^2 c^2}\Big)^{s_2}   \ll T^{-\delta},
\end{multline*}
for any $s_i$ with $\Re(s_i)=\ep$, $|\Im s_i|<T^\ep$, where $Z$ was slightly redefined. Absorbing the factors raised to the power $s_i$ into the existing weight functions, forgoing cancellation in the $k$ and $l$ sums, and treating only the negative sign case of $\sum_\pm$ as the positive case is similar, we are reduced to proving
\begin{align}
 \label{red1} \frac{N}{T^2}  \sum_{m, c \ge 1} \sums_{h \bmod c} e\Big(\frac{h(m + r)}{c}\Big)   \frac{  \tau(m,2T) }{ c^3}   Z\Big(\frac{\sqrt{Nm} }{cT}\Big)  \Psi_2\Big(\frac{m}{M}\Big)\ll T^{-\delta}
\end{align}
for any $r<T^\ep$. 
\subsection{Second application of Voronoi summation}

We now write the left hand side of (\ref{red1}) as
\begin{align*}
\frac{NM}{T^2}  \sum_{m, c \ge 1} \sums_{h \bmod c} e\Big(\frac{h(m + r)}{c}\Big)   \frac{  \tau(m,2T) }{ m c^3} \frac{m}{M} Z\Big(\frac{  \sqrt{Nm} }{cT}\Big)  \Psi_2\Big(\frac{m}{M}\Big)
\end{align*}
and apply Voronoi summation to the $m$-sum using Lemma \ref{vorsum2} with $\Phi(x)=x\Psi_2(x) Z(\frac{  \sqrt{xNM} }{cT})$. The main term $\mathcal{M}$ given by (\ref{vor2md}) is $O(T^{-A})$ for any $A>0$ by repeatedly integrating by parts. Thus it remains to consider the dual sum 
\begin{align}
\label{d2} \mathcal{D}=\frac{NM}{T^2} \sum_\pm \sum_{b, c \ge 1}  \frac{  \tau(b,2T) S(r, \pm b,c)}{ c^4}  \check{\Phi}^\pm\Big(\frac{Mb}{c^2}, T\Big).
\end{align}
In the definition of $\check{\Phi}^\pm (\frac{Mb}{c^2}, T)$ given in (\ref{phi-trans2}), we can restrict to $|s|<T^\ep$ by (\ref{phimell}) and (\ref{g2bound}). Now we can move the line of integration far to the right if $\frac{Mb}{c^2}>T^2$ and far to the left if $\frac{Mb}{c^2}<T^2$. Note that this does not cross any poles of $G_2^\pm(s)$ because of the restriction on $s$, and with (\ref{g2bound}) we have shown that we may restrict to
\begin{align*}
T^{2-\ep} < \frac{Mb}{c^2} <T^{2+\ep},
\end{align*}
or equivalently
\begin{align}
\label{crang} T^{-\ep} < \frac{b}{c^2} <T^{+\ep}
\end{align}
by (\ref{ranges0}). Thus the maximum size of $b$ is $T^{2+\ep}$, since $c<T^{1+\ep}$ by (\ref{ranges}). We now simplify the sum (\ref{d2}) as we did in the previous section. Introducing a smooth bump function $\Phi(\frac{4\pi \sqrt{rb}}{c})$ to make restriction (\ref{crang}), introducing a smooth bump function $\Psi(\frac{b}{B})$ to restrict to a dyadic interval $b\asymp B$, using (\ref{g2asymp2}) with $\Re(s)=\ep$ and $|\Im(s)|<T^\ep$, replacing a factor of $c^2$ with $N$ and a factor of $c$ with $b^{\half}$ by redefining the existing weight functions, replacing $M$ with $T^2$ at the front of the sum, and dropping the $Z$ function by absorbing it into the existing weight functions (by first separating variables using the Mellin transform), we are reduced to proving
\begin{align}
\label{vorend} \frac{1}{T} \sum_{b, c \ge 1}  \frac{  \tau(b,2T) S(r,  b,c)}{  b^{\half} c }  \Phi\Big(\frac{4\pi \sqrt{rb}}{c}\Big) \Psi\Big(\frac{b}{B}\Big) \Big(1+O(T^{-1+\ep})\Big) \ll T^{-\delta},
\end{align}
for any $r<T^\ep$ and $B<T^{2+\ep}$. We remind the reader that the factor of $\frac1T$ in front of the sum comes from (\ref{g2asymp2}). The total contribution of the error term is $O(T^{-\half+\ep})$, by using Weil's bound for the Kloosterman sum. Thus we are left to deal with the main term.

\subsection{Kuznetsov's formula and subconvexity}

The next step is to apply Kuznetsov's formula to the sum over $c$, in order to convert the sum of Kloosterman sums to a sum over automorphic forms. Using \cite[Lemma 3.5]{butkha} but written in terms of Hecke eigenvalues, we have
\begin{align}
\label{kuzbac} \sum_{c \ge 1}  \frac{  S(r,  b,c)}{   c }  \Phi \Big(\frac{4\pi \sqrt{rb}}{c} \Big) = \sum_{j\ge 1} \hat{\Phi}(t_j) \frac{4\pi|\rho_j(1)|^2}{\cosh(\pi t_j)} \lambda_j(r)\lambda_j(b) + \ldots,
\end{align}
where the sum is over an orthonormal basis of Hecke Maass cusp forms $\{u_j: j\ge1\}$ for $\Gamma$, with Laplacian eigenvalue $\frac{1}{4}+t_j^2$ (recall that there are no exceptional eigenvalues), $\hat{\Phi}(t)$ is the transform given in \cite[section 3.4]{butkha}, $\rho_j(1)$ is the first Fourier coefficient of $u_j$, which satisfies $|\rho_j(1)|^2\ll \cosh(\pi t_j)(1+|t_j|)^2$ by the standard (Rankin-Selberg) bound, and the ellipsis denotes the contribution of the holomorphic cusp forms and the Eisenstein series, which are fully written out in \cite[Lemma 3.5]{butkha} and whose treatment is similar. All we need to know about the transform is that for $t\in \mathbb{R}$, we have $\hat{\Phi}(t) \ll T^{-A}$ for any $A>0$ unless $|t|<T^\ep$, in which case $|\hat{\Phi}(t)| \ll T^{\ep}$.  This is shown in \cite[Lemma 3.6]{butkha}. Since there are only $O(T^\ep)$ forms $u_j$ to consider, we do not need cancellation from the $j$-sum.

Combining (\ref{vorend}) and (\ref{kuzbac}), we are reduced to proving (we write out the treatment for only the Maass forms) that
\begin{align}
\label{end1} \frac{1}{T} \sum_{b \ge 1}  \frac{  \tau(b,2T)  \lambda_j(b) }{  b^{\half}  }  \Psi\Big(\frac{b}{B}\Big) \ll T^{-\delta}
\end{align}
where $\lambda_j(n)$ are the Hecke eigenvalues of any $u_j$ with $|t_j|<T^\ep$, and $B<T^{2+\ep}$. Note that
\begin{align*}
L(\thalf + s+2iT, u_j) L(\thalf + s-2iT, u_j)= \zeta(1+2s) \sum_{b \ge 1}  \frac{  \tau(b,2T)  \lambda_j(n) }{  b^{\half+s}  }
\end{align*}
for $s>\half$. Thus by Mellin inversion, we have
\begin{align*}
\frac{1}{T} \sum_{b \ge 1}  \frac{  \tau(b,2T)  \lambda_j(b) }{  b^{\half}  }  \Psi\Big(\frac{b}{B}\Big) = \frac{1}{T} \frac{1}{2\pi i} \int_{(\sigma)} \frac{ L(\thalf + s+2iT, u_j) L(\thalf + s-2iT, u_j)}{ \zeta(1+2s)} B^{s} \tilde{\Psi}(s) ds
\end{align*}
for any $\sigma>\half$. We restrict the integral to $|\Im(s)|<T^\ep$ using the bound $|\tilde{\Psi}(s)|\ll T^\ep (\frac{T^\ep}{1+|s|})^{A}$ for any $A>0$ which follows by integration by parts. We then move the line of integration to $\sigma= 0$, where $|B^s| \ll 1$, $\zeta(1+2s)^{-1} \ll (1+|s|)^\ep$ as already noted earlier, and 
\begin{align*}
L(\thalf + s+2iT, u_j) L(\thalf + s-2iT, u_j) \ll T^{1-\delta}
\end{align*}
by any hybrid subconvexity bound (such as that of Jutila \cite{jut3}). This proves the bound required in (\ref{end1}).

\bibliographystyle{amsplain}

\bibliography{Eis}

\end{document}